\documentclass[12pt,leqno]{article}

\usepackage[hypertex,backref]{hyperref}
\usepackage{amsmath}
\usepackage{amssymb}

\usepackage[american]{babel}

\usepackage{eucal,euler,times}

\usepackage{amsthm}
\usepackage{enumerate}

\usepackage{xy}
\xyoption{all}

\swapnumbers

\theoremstyle{plain}
    \newtheorem{theorem}{Theorem}[section]
    \newtheorem{lemma}[theorem]{Lemma}
    
    \newtheorem{proposition}[theorem]{Proposition}

\theoremstyle{definition}
    \newtheorem{definition}[theorem]{Definition}
    \newtheorem{example}[theorem]{Example}

     \newtheorem{remark}[theorem]{Remark}

\theoremstyle{remark}

\numberwithin{equation}{section}

\def\R{{\mathbb R}}
\def\Z{{\mathbb Z}}

\def\C{{\mathbb C}}

\def\S{{S}}

\def\K{\mathbb{K}}

\def\Bounded{\mathcal{B}}
\DeclareMathOperator{\Spin}{Spin}

\DeclareMathOperator{\ev}{ev}
\DeclareMathOperator{\geom}{geom}
\DeclareMathOperator{\opp}{opp}
\DeclareMathOperator{\Pt}{pt}

\DeclareMathOperator{\odd}{odd}
\DeclareMathOperator{\Cliff}{Cliff}
\DeclareMathOperator{\rank}{rank}

\DeclareMathOperator{\End}{End}

\begin{document}

\title{On the Equivalence of  Geometric and Analytic K-Homology}
\author{Paul Baum, Nigel Higson, and Thomas Schick}
\date{\today}

\maketitle

\begin{abstract}
We give a proof that the geometric $K$-homology theory  for finite $CW$-complexes defined by Baum and Douglas
is isomorphic to Kasparov's $K$-homology.  The proof is a
simplification of more elaborate arguments which deal with the
geometric formulation of \emph{equivariant} $K$-homology theory.
\end{abstract}

\section{Introduction}

$K$-homology theory, the  homology theory which is dual to
Atiyah-Hirzebruch $K$-theory, may be defined abstractly using the Bott
spectrum and standard contructions in homotopy theory.
Atiyah \cite{atiyah69} pointed out the relevance to index theory of a concrete
definition of $K$-homology. Following his suggestions, detailed
analytic definitions of $K$-homology were provided by Brown, Douglas
and Fillmore \cite{BDF77} and by Kasparov \cite{kasparov75}, and these works are now
foundational papers in operator $K$-theory.  At about the same
time, Baum and Douglas  \cite{baum-douglas82} introduced a geometric definition of $K$-homology (using manifolds, bordisms, and so on) in connection with work on the
Riemann-Roch problem \cite{MR54:317,MR82c:14021}.  Baum and Douglas  defined a very simple
and natural map from their geometric theory to analytic $K$-homology, and
this map turns out to be an isomorphism.  The combined efforts of
various mathematicians in the early 1980's produced a proof of this,
but a detailed account of the matter was never published. This is
despite the fact that over the years the isomorphism has grown in
importance, thanks to its connection with the Baum-Connes conjecture
\cite{baum-connes-higson94}.  The purpose of this note is to present, after a twenty five year
gap, a detailed proof of the isomorphism from geometric
$K$-homology to analytic $K$-homology.  (See \cite{jakob1,jakob2} for a related approach to the problem of defining homology theories dual to multiplicative cohomology theories like $K$-theory.)
The proof is a spin-off from
our work on equivariant $K$-homology theory, which will be reported
upon in a future paper, where  we shall prove that for a discrete, countable group $G$,  geometric equivariant $K$-homology is isomorphic to analytic equivariant $K$-homology on the category of proper, finite $G$-$CW$-complexes.

With admiration and affection we dedicate this paper to Robert MacPherson.  A conversation between the first-named author and Bob MacPherson at IHES in 1978 was crucial to the eventual formulation of geometric $K$-homology.

\section{Review of Analytic K-Homology} 
\label{kasp-sec}

In this section we shall review Kasparov's definition of analytic $K$-homology,
and list those facts about it  that we shall need in
the sequel.  For further details the reader is referred to the
monograph \cite{MR2002c:58036} on the subject.

Throughout this section we shall be working with locally compact, second countable topological spaces.  If $Z$ is such a
space then we shall denote by $C_0(Z)$ the (separable) $C^*$-algebra of
continuous, complex-valued functions on $Z$ which vanish at infinity.

If $X$ and $Y$ are operators on a Hilbert space, then the notation $X \sim
Y$ will signify the equality of $X$ and $Y$ modulo the compact operators.

\begin{definition}
\label{fredmoddef}   
Let $A$ be a separable $C^*$-algebra.  An (ungraded) \emph{Fredholm
module} over $A$ is given by the following data:
\begin{enumerate}[\rm (a)] 
  
\item a separable Hilbert space $H$,

\item  a  representation $\rho\colon A\to
\Bounded(H)$ of $A$ as bounded operators on $H$,  and

\item  an operator $F$ on $H$ such that for all $a\in A$,
$$
(F^2- 1)\rho(a)\sim 0,\quad (F-F^*)\rho(a)\sim 0,
\quad F\rho(a)\sim \rho(a)F.
$$ 
\end{enumerate} 
The representation $\rho$ is not required to be non-degenerate in any
way.  In fact $\rho$, and even the Hilbert space $H$, are allowed to be
zero.
\end{definition}

Roughly speaking, Kasparov's $K$-homology groups are assembled from
homotopy classes of Fredholm modules over $A=C_0(Z)$.  However it is
necessary to equip these Fredholm modules   with a modest
amount of extra structure.

\begin{definition}
Let $p\in \{0,1,2,\dots \}$ and let $A$ be a separable $C^*$-algebra.
A \emph{$p$-graded}\footnote{The term `$p$-multigraded' is used in
\cite{MR2002c:58036}.}  Fredholm module is a Fredholm module $(H,\rho, F)$, as
above, with the following additional structure:

\begin{enumerate}[\rm (a)]

\item The Hilbert space $H$ is equipped with a $\Z/2$-grading
$H=H^+\oplus H^-$ in such a way that  for each $a\in A$, the operator
$\rho(a)$ is even-graded, while the operator $F$ is odd-graded.

\item  There are  odd-graded operators
$\varepsilon_1,\dots , \varepsilon_p$ on $H$ such that
$$ 
        \varepsilon_j = -\varepsilon_j^*,\quad \varepsilon_j^2 =
        -1,\quad \varepsilon_i\varepsilon_j+\varepsilon_j\varepsilon_i
        = 0\quad (i\ne j),
$$ 
and such that $F$ and each $\rho (a)$ commute with each $\varepsilon_j$.
 \end{enumerate} 
\end{definition}

Of course, if $p=0$ then part (b) of the definition does not apply. 

\begin{definition}
Let $(H,\rho,F)$ and $(H',\rho', F') $ be $p$-graded Fredholm modules
over $A$. A \emph{unitary equivalence} between them is a
grading-degree zero unitary isomorphism $U\colon H\to H'$ which
intertwines the representations $\rho$ and $\rho'$, the operators $F$
and $F'$, and the grading operators $\varepsilon_j$ and
$\varepsilon_j'$.
\end{definition}

\begin{definition}
Suppose that $(H,\rho,F_t)$ is a family of $p$-graded Fredholm modules
parameterized by $t\in [0,1]$, in which the representation $\rho$,
the Hilbert space $H$ and its grading structures remain constant but
the operator $F_t$ varies with $t$.  If the function $t\mapsto F_t$ is
norm continuous, then we say that the family defines an \emph{operator
homotopy} between the $p$-graded Fredholm modules $(\rho,H,F_0)$ and
$(\rho,H,F_1)$, and that the two Fredholm modules  are \emph{operator homotopic}.
\end{definition}

There is a natural notion of \emph{direct sum} for Fredholm modules: one
takes the direct sum of the Hilbert spaces, of the representations, and of the
operators $F$.  The \emph{zero module}  has zero Hilbert space, zero
representation, and zero operator.

Now we can give Kasparov's definition of $K$-homology.

\begin{definition}
Let $p\in \{0,1,2,\dots \}$ and let $A$ be a separable $C^*$-algebra.
The \emph{Kasparov $K$-homology group} $K^{-p}(A)$ is the abelian
group with one generator $[x]$ for each unitary equivalence class of
$p$-graded Fredholm modules over $A$ and with the following relations:

\begin{enumerate}[\rm (a)]

\item  if $x_0$ and $x_1$
are operator homotopic $p$-graded Fredholm modules then $[x_0]=[x_1]$ in
$K^{-p}(A)$, and 

\item 
if $x_0$ and $x_1$ are any two $p$-graded Fredholm modules then
$[x_0\oplus x_1] = [x_0] + [x_1]$ in $K^{-p}(A)$.
\end{enumerate}  
\end{definition}

\begin{definition}
A $p$-graded Fredholm module is said to be \emph{degenerate} if the equivalences modulo compact operators listed in item (c) of Definition~\ref{fredmoddef} are actually equalities. 
\end{definition}

It is easy to see that a degenerate $p$-graded Fredholm module determines the zero element of $K^{-p}(A)$. 

\begin{lemma}
\label{triv-mod-lemma}
Let $(H,\rho, F)$ be a $p$-graded Fredholm module.  Assume that there exists a self-adjoint, odd-graded involution $E\colon H\to H$ which commutes with the action of $A$ and with the multigrading operators $\varepsilon _j$, and which anticommutes with $F$.  Then the Fredholm module $(H,\rho, F)$ represents the zero element of $K^{-p}(A)$.
\end{lemma}

\begin{proof}
The path $F_t = \cos (t) F + \sin (t) E$ gives an operator homotopy from $F$ to the degenerate operator $E$.
\end{proof}

It follows from the lemma  that the additive inverse of the $K$-homology class represented by $(H,\rho , F)$ is the class of $(H^{\opp}, \rho, -F)$, where $H^{\opp}$ denotes $H$ with the grading reversed.  This is because the involution $\left ( \begin{smallmatrix} 0 & 1 \\ 1 & 0 \end{smallmatrix} \right )$ on $H\oplus H^{\opp}$  satisfies the hypotheses of the lemma, applied to the Fredholm module $(H\oplus H^{\opp},\rho\oplus \rho,F\oplus- F)$. 
It follows that every class in $K^{-p}(A)$ is represented by a single Fredholm module, and that two modules represent the same class if and only if, up to isomorphism, they become operator homotopic after adding degenerate modules.

If $(H,\rho, F)$ is a $p$-graded Fredholm module $A$, then we may construct
from it a $(p+2)$-graded Fredholm module $(H',\rho',F')$ over $A$ by
means of the formulas 
$$
        H' = H\oplus H^{\opp}, \quad \rho' = \rho\oplus \rho, \quad F'
        = F \oplus F,
$$
 along with the
grading operators
$$
 \varepsilon _j = \varepsilon _j\oplus \varepsilon_j\quad (j=1,\dots
 p),
\quad 
\varepsilon _{p+1} = \begin{pmatrix} 0 & I \\ -I& 0 \end{pmatrix}\quad
\text{and}\quad \varepsilon _{p+2} = \begin{pmatrix} 0 & iI \\ iI & 0
\end{pmatrix}.
$$
\begin{definition} 
\label{periodicity-def}
The \emph{formal periodicity map} 
$$
     K^{-p}(A) \xrightarrow{} K^{-(p+2)}(A)
$$
is the homomorphism of Kasparov  groups induced from this construction.
\end{definition} 

The periodicity map can be reversed by compressing a $(p+2)$ graded
Fredholm module to the $+1$ eigenspace of the involution $-i
\varepsilon _{p+1}\varepsilon _{p+2}$.  We obtain an isomorphism 
$$
        K^{-p}(A) \cong K^{-(p+2)}(A).
$$
As a result there are really only two genuinely distinct $K$-homology
groups, $K^{\ev}$ and $K^{\odd}$, as follows:

\begin{definition}
Let us denote by $K^{\ev}(A)$ and  $K^{\odd}(A)$  the groups
$K^0(A)$ and $K^{-1}(A)$ respectively, or more canonically, the direct
limits 
$$ K^{\ev}(A) = \varinjlim _k K^{-{2k}}(A) \quad \text{and}
\quad K^{\odd}(A) =  \varinjlim _k K^{-(1+2k)}(A)$$
under the above periodicity maps. \end{definition}

\begin{definition}
If $Z$ is a second countable,\footnote{This assumption is required at several points in Kasparov's theory, which is designed for \emph{separable} $C^*$-algebras.}  locally compact space, and if  $A= C_0(Z)$, then we shall write $K_{p}(Z) $ in place of $K^{-p}(A)$.
These are the \emph{Kasparov $K$-homology groups} of the space $Z$. 
If $(X,Y)$ is
a second countable, locally compact pair, and if $Z$ is the difference $X\setminus Y$, then we define relative $K$-homology groups by
$$
        K_{p} (X,Y) = K_{-p} (Z).
$$
We shall define periodic groups $K_{\ev/\odd}(X,Y)$ similarly. 
\end{definition}

Kasparov's main theorem concerning these objects is then as follows:

\begin{theorem}
There are natural transformations
$$
        \partial \colon K_{p} (X,Y) \xrightarrow{} K_{p-1} (Y)
$$
(connecting homomorphisms) which are compatible with the formal periodicity isomorphisms and which give Kasparov $K$-homology the
structure of a $\Z/2$-graded homology theory on the category of
compact metrizable pairs $(X,Y)$.  On the subcategory of finite
$CW$-complexes Kasparov $K$-homology is isomorphic to topological
$K$-homology --- the  homology theory associated to the
Bott spectrum. \qed
\end{theorem}

\section{Dirac-Type Operators}
\label{dirac-op-sec}

We continue to follow the monograph \cite{MR2002c:58036}. 

\begin{definition}
\label{diracbundle} 
Let $M$ be a smooth, second countable finite dimensional manifold (possibly
with non-empty boundary) and let  $V
$ be a smooth, Euclidean vector bundle over $M$. A \emph{$p$-graded Dirac
  structure} on 
$V $ is a smooth, $\Z/2$-graded, Hermitian vector bundle $\S$ over $M$ together
with the following data:
\begin{enumerate}[\rm (a)]
\item An $\R$-linear morphism of vector bundles $$ V  \to
\End(\S) $$  which associates to each vector $v\in V _x$ a
skew-adjoint, odd-graded endomorphism $u\mapsto v \cdot u$ of $\S_x$ in such
a way that 
$$
v\cdot v\cdot u = -\|v \|^2 u.
$$
\item
A family of skew-adjoint, odd-graded endomorphisms 
$\varepsilon_1,\dots , \varepsilon_p$ of $\S$ such that 
$$ 
        \varepsilon_j = -\varepsilon_j^*,\quad \varepsilon_j^2 =
        -1,\quad \varepsilon_i\varepsilon_j+\varepsilon_j\varepsilon_i
        = 0\quad (i\ne j),
$$ 
and such that each $\varepsilon _j$ commutes with each operator
$u\mapsto v \cdot u$.
\end{enumerate}
Usually $M$ will be a Riemannian manifold and we will take $V = T M$.  In this case we shall call $\S$  a \emph{$p$-graded Dirac bundle on $M$}.
\end{definition}

 \begin{definition}
\label{dirac-op-def}
Let $M$ be a Riemannian manifold which is equipped with a $p$-graded Dirac structure, with Dirac bundle $\S$. 
We shall call an  odd-graded, symmetric, order one linear partial differential operator $D$ acting on the sections of  $\S$ a \emph{Dirac operator} if it  commutes with the operators $\varepsilon _j$, and if  
$$
[ D, f] u = \operatorname{grad} f\cdot u,
$$
for every smooth function $f$ on $M$ and every section $u$ of $\S$. \end{definition}

Every Dirac bundle on a Riemannian manifold  admits a Dirac operator, and the difference of two Dirac operators on a single Dirac bundle $\S$ is an endomorphism of $\S$.

A $p$-graded Dirac operator $D$ on a Riemannian manifold $M$ without boundary
defines in a natural way a class $[D]\in K_p(M)$.  The general construction is a little involved, and we refer the reader to \cite{MR2002c:58036} for details, but when $M$ is closed there is a very simple description of $[D]$:

\begin{theorem} 
\label{fmod1}
Let $M$ be a closed (i.e.~compact without boundary) Riemannian manifold and let $D$ be a  Dirac operator on a $p$-graded Dirac bundle  $\S$. Let $H=L^2(M,S)$ be the Hilbert  space of square-integrable sections of $\S$, and let $\rho$ be the representation of
$C(M)$ on $H$ by pointwise multiplication operators. Let 
$$
        F = D ( I + D^2 )^{-\frac 12} .
$$ 
The triple $(\rho, H , F)$ is a $p$-graded Fredholm module for $A=
C(M)$. \qed 
\end{theorem}

To describe further properties of the classes $[D]$ we need to introduce the following boundary operation on Dirac bundles:

\begin{definition}
\label{boundary-S-def}
Let $\S$ be a $p$-graded Dirac bundle on a  Riemannian
manifold $\overline{M}$ with boundary $\partial M$. If $e_1$ denotes the outward 
pointing unit normal  vector field on the boundary manifold $\partial
M$ then the formula
$$
        X\colon u \mapsto (-1)^{\partial u}e_ 1\cdot \varepsilon_1 u 
$$ 
defines an automorphism of the restriction of $\S$ to $\partial M$
which is even, self-adjoint, and satisfies $X^2=1$. The operator $X$ commutes
with  multiplication $u\mapsto Y\cdot u$ by  tangent vectors  $Y$ orthogonal to $e_1$,
and also with the multigrading operators $\varepsilon_2, \dots ,
\varepsilon_{p}$.  The $+1$ eigenbundle for $X$ is a $(p-1)$-graded Dirac 
bundle\footnote{The multigrading operators are obtained from
$\varepsilon_2,\dots , \varepsilon_p$ by shifting indices downwards.} on
$\partial M$, which we shall call the \emph{boundary} of the Dirac
bundle $\S$.
\end{definition}

The following theorem summarizes facts proved in Chapters 10 and 11 of \cite{MR2002c:58036}.  

\begin{theorem}  
\label{bdiracdirac} 
To each Dirac operator $D$ on a $p$-graded Dirac bundle over a smooth manifold without boundary there is associated a class $[D]\in K_p (M)$ with the following properties: 
\begin{enumerate}[\rm (i)]
\item The class $[D]$ depends only on the Dirac bundle, not on the choice of the operator $D$.
\item 
If $M_1$ is an open subset of $M_2$, and if $D_1$ is a Dirac operator on $M_1$ obtained by restricting a Dirac operator $D_2$ on $M_2$, then $[D_2]$ maps to $[D_1]$ under the homomorphism $K_p(M_2) \to K_p (M_1)$.
\item
Let $M$ be the interior of a Riemannian manifold $\overline{M}$ with
boundary $\partial M$, and let  $\S$ be a $p$-graded Dirac bundle on
$\overline{M}$.  Let $D$ be a Dirac operator on $M$  associated to $\S$
and let $D_{\partial M}$ be a Dirac operator on $\partial M$
associated to the boundary of $\S$. The connecting homomorphism 
$$ \partial\colon K_p(M)\to K_{p-1}(\partial M) $$
in Kasparov $K$-homology  takes the
class $[D_M] $ to the class $[D_{\partial M}]$: 
$$ 
\partial [D_M] = [D_{\partial M}] \in K_{p-1}(\partial M) . 
$$
\end{enumerate}
\qed
\end{theorem}

We shall need one additional fact about Dirac operators which concerns the structure of operators  on fiber bundles.  Suppose that $M$ is a closed Riemannian manifold and that  $P$ is a principal bundle over $M$ whose structure group is a compact Lie group $G$.  Suppose that $N$ is a closed Riemannian manifold on which $G$ acts by isometries.  We can then form the manifold 
$ Z =P\times _G N$.  Its   tangent bundle $T Z$ 
fits into an exact sequence of vector bundles over $Z$,
$$
\xymatrix{
0 \ar[r] &V \ar[r] & T Z  \ar[r] &  \pi^* T M \ar[r] & 0},
$$
where $\pi$ denotes the projection mapping from $Z $ to $M$ and where $V$ denotes the ``vertical  tangent bundle'' $V = P\times_G T N$.
If we choose a splitting of the sequence then we obtain an isomorphism 
\begin{equation}
\label{dir-sum-eqn}
 T Z \cong V\oplus \pi^* T M  ,
\end{equation}
which equips $Z$ with a Riemannian metric. 

Now suppose that $\S_M$ is a $p$-graded Dirac bundle for $M$ and that $\S_N$ is a $0$-graded Dirac bundle for $N$.  Let us also suppose that there is an action of $G$ on $\S_N$ which is compatible with the action of $G$ on  $N$.   We can then form the bundle $\S_{ V}  = P\times _G \S_N$ over $Z$, and from it the graded tensor product 
$\S_Z =  \S_{ V}  \hat \otimes \pi^*   \S _M $.   Using the direct sum decomposition (\ref{dir-sum-eqn}) this becomes a $p$-graded Dirac bundle for $Z$, with the  tangent vector $v\oplus w\in V\oplus  \pi^* T M  $ acting as the operator $v\hat\otimes 1 + 1 \hat \otimes w
 $ on $ \S_{V}\hat \otimes \pi^*  \S _M  $. 

We can now form the class $[D_{Z}]\in K_{p} (Z)$ associated to a Dirac operator on the Dirac bundle $\S_{V}\hat \otimes \pi^*   \S _M  $, and using the projection mapping $\pi\colon Z\to M$ we obtain a class 
$$
\pi_*[D_{Z}] \in K_p (M).
$$
The following proposition relates $\pi_* [D_{Z}] $ to the class
$[D_M]$ of a Dirac operator for the Dirac bundle $\S_M$ on $M$.

\begin{proposition} 
\label{fiber-bundle-prop}
Assume 
that there exists a $G$-equivariant Dirac operator for the Dirac bundle $\S_N$ on $N$ whose kernel is the one-dimensional trivial representation of $G$, spanned by an even-graded section of $\S_N$. 
Then 
$$
\pi_* [D_{Z}] = [D_M] \in K_p(M).
$$
\end{proposition}

\begin{proof}
Let us consider first the special case in which the principal bundle $P$ is trivial: $P = G\times M$ (in this case we might as well take $G= \{ e\}$).  Then of course $Z = N \times M$. 
We can take the Dirac operator $D_{Z}$ to be 
$$
 D_{Z}= D_N \hat \otimes I + I \hat \otimes D_M,
 $$
where $D_N$ is a Dirac operator for the Dirac bundle $\S$ on $N$ with one-dimensional kernel, as in the statement of the proposition.  Now the Hilbert space on which $D_{Z}$ acts is the tensor product 
$$
L^2(N\times M , \S_N\hat \otimes \S_M) = L^2 (N, \S_N) \hat \otimes  L^2 (M,\S_M) .
 $$
If we split the first  factor, $L^2 (N,\S_N)$, as $\ker (D_{N})$ plus its orthogonal complement, then we obtain a corresponding direct sum decomposition of $L^2(N\times M , \S_N\hat \otimes \S_M)$.   The operator $F_{Z}$ formed from $D_{Z}$, as in Theorem~\ref{fmod1}, respects this direct sum decomposition, as does the action of $C(M)$.  We therefore obtain a decomposition of the Fredholm module representing $[D_{Z}]$ as a direct sum of two Fredholm modules.  The first acts on $ \ker (D_N) \otimes L^2 (M, \S_M)\cong L^2 (M,\S_M)$ and is isomorphic to the Fredholm module representing $[D_M]$.  The second represents the zero element of $K_p(M)$.  This follows from Lemma~\ref{triv-mod-lemma}, since if $T$ is the partial isometry part of $D_N$ in the polar decomposition, and if $\gamma$ is the grading operator on $L^2 (M,\S_M)$,  then the odd-graded involution 
$$
 E = i T \hat \otimes   1 
$$
on the Hilbert space $  \ker (D)^{\perp}\hat \otimes L^2 (M,\S_M) $ commutes
with the action of $C(M)$, and with the grading operators $\varepsilon _j$,
and anticommutes with $F_{Z}$ (using the conventions of \cite[A 3.3]{MR2002c:58036}.

The proof of the general case is similar.  To begin, the Hilbert space on which $D_{Z}$ acts is naturally isomorphic to the fixed point space 
$$
 \left [L^2 (N,\S_N) \hat \otimes L^2 (P, \pi^* \S_M) \right ] ^G.
$$
Denote by $\widetilde D_M$ a $G$-equivariant  linear partial differential operator on $P$, acting on sections of $\pi^* \S_M$, which is obtained
as follows.  Select a finite cover of $M$ by open sets $U_j$ over which the bundle $P$ is trivial, and fix isomorphisms to $G\times U_j$ over these open sets.  Use the isomorphisms to define operators $\widetilde D_j $ on $\pi^{-1}[U_j]\subseteq P$ which act as $D$ in  the $U_j$ direction and act as the identity in the $G$-direction.  Select also a smooth partition of unity $\{ \sigma_j^2\}$ which is subordinate to the cover.  Then define $\widetilde D$ by averaging the sum $ \sum  \sigma _j \widetilde D_j\sigma_j $ over the action of $G$.
 Having constructed $\widetilde D_M$, we obtain a Dirac operator for $\S_M\hat \otimes \S_V$ by the formula 
$$
D_{Z} = D_N \otimes I + I \hat \otimes  \widetilde D_M.
$$
From here the argument used in the special case may be applied \emph{verbatim}.\end{proof}

\begin{remark}
By using some machinery the preceding result can be conceptualized and generalized  as follows.  If $G$ is a compact group and $A$ is a $C^*$-algebra equipped with an action of $G$ (for example $A=C(N)$), then there is a natural notion of $G$-equivariant Fredholm module, from which we may define equivariant $K$-homology groups $K^{-p}_G (A)$.  In the commutative case these give equivariant groups $K^G _{p} (N)$.  Now if $P$ is a principal $G$-bundle over $M$, as above, then by elaborating on the construction of the Kasparov product (which we shall not actually use anywhere in this paper) we obtain a pairing 
$$
\xymatrix{
 K^G_0(N)\otimes K_{p}(M) \ar[r]^-{\mu_P}& K_{p} (Z)}
$$
One can compute that the class $[D_M]\otimes [D]$ is mapped to $[D_{Z}]$.
Next, the map which collapses $N$ to a point induces a homomorphism from $K^G_0(N) $ to the coefficient group $K_0^G(\text{pt})$, which is the representation ring of $G$.  From a representation of $G$ and the principal bundle $P$ we obtain by induction a vector bundle on the space $M$.  We therefore  obtain a map 
$$
\xymatrix{
\varepsilon \colon K_0^G(N) \ar[r]&K^0(M)}.
$$
Finally, the group $K_p(M)$ is a module over the ring $K^0(M)$ by the cap product between homology and cohomology.  We obtain a diagram
$$
\xymatrix{
 K^G_0(N)\otimes K_{p}(M)\ar[r]^-{\mu_P}\ar[d]_{  \varepsilon\otimes 1} & K_{p} (Z)\ar[d]^{\pi_*}\ \\
K^{0} (M) \otimes K_{p} (M)\ar[r]_-{\cap} & K_{p}(M). }  
$$
Proposition~\ref{fiber-bundle-prop}  follows from the assertion that this diagram commutes (in the special case where the collapse map sends $[D]$ to $1\in R(G)$).  The commutativity of the diagram is a simple exercise with the Kasparov product, but it is beyond the scope of the present article. 
\end{remark}

We conclude this section by introducing a specific Dirac  operator to which we shall apply Proposition~\ref{fiber-bundle-prop}.   In order to fix notation we begin with the following definition:

\begin{definition}
Let $V$ be a Euclidean vector space. The \emph{complex Clifford algebra} for $V$
is the universal complex $*$-algebra $\Cliff (V)$ equipped with an $\R$-linear inclusion of $V$, and subject to the relations 
$
v^2 = -\| v\|^2 \cdot 1$ for $v\in V$.  If $\{ e_1,\dots, e_n\}$ is an orthonormal basis for $V$, then the algebra $\Cliff (V)$ is linearly spanned by the $2^n$
monomials $e_{j_1}\cdots e_{j_k}$, where $j_1<\dots <j_k$ and $0\le k \le n$.
We introduce an inner product on $\Cliff (V)$ by deeming these monomials to be
orthonormal. 
\end{definition}

The algebra $\Cliff (V)$ is $\Z/2$-graded: the monomial $e_{j_1}\cdots e_{j_k}$ is even or odd-graded, according as $k$ is even or odd.

\begin{definition}
Let $N$ be an even-dimensional, Riemannian manifold and let  $\Cliff(TN)$ be the complex vector bundle on $N$ whose fibers are the complexified Clifford algebras of the fibers of the tangent bundle of $N$.  The 
bundle $\Cliff (TN)$ has a natural $0$-graded Dirac bundle structure (tangent vectors act by Clifford multiplication on the left).
\end{definition} 

If  $N$ is oriented, and if $\{e_1,\dots, e_n\}$ is a local, oriented, orthonormal frame, then the operator of \emph{right}-multiplicaton by the product 
$$
 \sigma = i^{\frac n2 } e_1\cdots e_n
$$
is an even-graded, self-adjoint involution of the bundle $\Cliff (TN)$ which commutes with the Dirac bundle structure.  

\begin{definition}
Denote by $\Cliff_{\frac 12} (TN)$ the $+1$-eigenbundle of the involution $\sigma$.  This is a $0$-graded Dirac bundle in its own right. 
\end{definition}

We wish to compute the index of a Dirac operator associated to this Dirac bundle, at least in the case of a sphere $N=S^n$.   To do so, we use the standard isomorphism between 
$
\Cliff (TN)$ and the complexified exterior algebra bundle $\bigwedge^*_\C T^*N\cong \bigwedge^*_\C T N$, which associates to the Clifford monomial $e_{j_1}\cdots e_{j_k}$ the differential form $e_{j_1} \wedge \dots \wedge e_{j_k}$.  Under  this correspondence, the operator $D=d+d^*$ on forms becomes a Dirac operator for the Dirac bundle $\Cliff (TN)$.  So the kernel of $D$ is the space of harmonic forms on $N$.  Using the fact that the involution $\sigma$ exchanges the $0$ and $n$-forms on $N$ we obtain the following result. 

\begin{proposition}
\label{index-calc-prop}
Let $N$ be an even-dimensional, round sphere (oriented as the boundary of the ball).  There is a Dirac operator for $\Cliff _{\frac 12 } (TN)$ which is equivariant for the natural action of the special orthogonal group, and whose kernel is the one-dimensional trivial representation, and is generated by an even-graded section of $\Cliff_{\frac 12}(TN)$. \qed
\end{proposition}

\begin{remark}
For general oriented Riemannian manifolds $N$, the index of the Dirac operator for $\Cliff _{\frac 12 } (TN)$  is the average of the Euler characteristic and the signature.
Indeed the direct sum of $\Cliff_{\frac 12} (TN)$ with the opposite of the bundle complementary to $\Cliff _{\frac 12}(TN)$ in $\Cliff (TN)$ is the Dirac bundle associated to the signature operator of Atiyah and Singer. 
\end{remark}

\section{Spin$^c$-Structures}
\label{spinc-sec} 
We shall define $\Spin^c$-structures using the notion of Dirac bundle that was introduced in the last section.  

\begin{definition}
Denote by $\C_n$ the  complex Clifford algebra for $\R^n$,  
generated by the standard basis elements $e_1,\dots , e_n$ of $\R^n$.
\end{definition}

Let  $M$ be a smooth manifold and let $V$ be a  rank $p$
Euclidean vector bundle over $M$. If $e_1,\dots ,
e_n$ is a local orthonormal frame for $V$, defined over an open set
$U\subseteq M$, then the trivial bundle $U\times \C_n$ over $U$ with
fiber $\C_n$ may be given the structure of an $p$-graded Dirac bundle for $V\vert _U$:
Clifford multiplication by an element $e_j$ of the frame is
\emph{left} multiplication by the $j$th generator of $\C_p$, and the
$p$-multigrading operators $\varepsilon_1,\dots , \varepsilon_p$ for the
bundle are \emph{right} multiplication by the same
generators.

\begin{definition} 
\label{spinordef}  
Let $M$ be a smooth manifold and let $V$ be a  $p$-dimensional
Euclidean vector bundle over $M$. A \emph{complex spinor bundle} for $V$ is a 
$p$-multigraded Dirac bundle $\S_V$ which is locally isomorphic to the trivial
bundle with fiber $\C_p$, the Clifford multiplication being determined from
some local orthonormal frame, as above.  We shall call a bundle $V$ equipped with a complex spinor bundle a \emph{$\Spin^c$-vector bundle}.   If $M$ is a smooth manifold (possibly with boundary) then by a \emph{$\Spin^c$-structure} on $M$ we shall mean a pair consisting of a Riemannian metric on $M$ and a complex spinor bundle $\S_M$ for $T M$.
\end{definition}

\begin{remark}
A spinor bundle determines an orientation of $V$, as follows. 
If $\{ f_1,\dots, f_p\}$ is a local orthonormal frame for $V$, then the endomorphism of the spinor bundle $\S_V$ determined by the formula 
$$ 
        u\mapsto (-1)^{p+(p-1)\partial
        u} f_1\dotsm f_p\cdot \varepsilon_p\dotsm \varepsilon _1 u 
$$ 
is plus or minus the identity (here $\partial u$ is the $\Z/2$-grading degree of the section $u$). If the endomorphism is $+I$ then we deem  the frame to be oriented; if it is $-I$ then we deem it to be oppositely oriented. 
\end{remark}

\begin{example}
Let $V_1$ and $V_2$ be Euclidean vector bundles on $M$ equipped with spinor bundles $\S_1$ and $\S_2$. Using the well-known Clifford algebra isomorphism $\C_{p_1}\hat \otimes \C_{p_2} \cong \C_{p_1+p_2}$ the graded tensor product $S_1\hat\otimes S_2$ becomes a spinor bundle for $V_1\oplus V_2$.   It defines  the \emph{direct sum}  $\Spin^c$-structure on $V_1\oplus V_2$. 
\end{example}

\begin{remark}
The definition of $\Spin^c$-structure can be rephrased in the language of principal bundles, as follows.  The group $\Spin(n)$ is the closed subgroup of the unitary group of $\C_n$ whose Lie algebra is the $\R$-linear span of the elements $e_ie_j$, for $i\ne j$.   The group $\Spin^c (n)$ is the closed subgroup of the unitary group of $\C_n$ which is generated by $\Spin(n)$ and the 
complex numbers of modulus one.  The group $\Spin^c(n)$ acts by inner automorphisms on the $\R$-linear subspace of $\C_n$ spanned by the elements $e_j$, and in this way we obtain a homomorphism from $\Spin^c(n)$ into $GL(n,\R)$ (in fact into $O(n)$). 
Now if $M$ is a smooth manifold, and if $  P$ is a reduction to $\Spin^c (n)$ of the principal bundle of tangent frames, then the reduction determines a Riemannian metric on $M$, and the bundle
$$
 \S = P \times _{\Spin^c (n)} \C_n
$$
is a spinor bundle on $M$ (here $\Spin^c (n)$ acts on $\C_n$ by left multiplication).  Thus $P$ determines a $\Spin^c$-structure.  Conversely, every $\Spin^c$-structure arises in this way (up to isomorphism).
\end{remark}

\begin{definition} 
Let $M^n$ be a smooth manifold, without boundary, equipped with a $\Spin^c$-structure.  We shall denote by $[M]\in K_n(M)$ the $K$-homology class of any Dirac operator on $\S$.  This is the $K$-homology \emph{fundamental class} of the $\Spin^c$-manifold $M$.
\end{definition}

If $\overline{M}$ is a smooth manifold with boundary then of course a Riemannian metric on $\overline{M}$ restricts to one on the interior $M$, and also to one on the boundary $\partial M$.  A spinor bundle $S$ for  $\overline M$ restricts to a spinor bundle on $M$, and the boundary of $S$, as described in Definition~\ref{boundary-S-def}, is a spinor bundle for $\partial M$.   The following result is a consequence of Theorem~\ref{bdiracdirac}.

\begin{theorem}  
\label{bdiracdirac2} 
If $M$ is the interior of an $n$-dimensional
$\Spin^c$-manifold with boundary, and if we equip the boundary manifold
$\partial M$ with the induced $\Spin^c$-structure, then the $K$-homology
boundary map $$ \partial\colon K_n(M)\to K_{n-1}(\partial M) $$  takes the
fundamental class of $M$ to the fundamental class of $\partial M$: $$ \partial
[M] = [\partial M] \in K_{n-1}(\partial M) . $$ 
\qed
\end{theorem}

\begin{definition}
\label{opp-spin-c-def}
Let $M$ be a smooth manifold equipped with a $\Spin^c$-structure.  The \emph{opposite} $\Spin^c$-structure is defined by changing the action of the multigrading operator $\varepsilon _1$ by a sign. 
\end{definition}

\begin{definition}  
Let $M$ be a smooth manifold.  Two $\Spin^c$-structures on $M$ are 
\emph{concordant} if there is a $\Spin^c$-structure on  $[0,1] \times M$ for which the induced $\Spin^c$-structure on the boundary $M \cup M$ is one of the given $\Spin^c$ structures on one copy of $M$, and the opposite of the other given structure on the other copy of $M$.
\end{definition}

In Chapter 11 of \cite{MR2002c:58036}, the following result is proved.
\begin{theorem}  
Concordant $\Spin^c$-structures on $M$ determine the same fundamental class in $K$-homology.\quad
\qed 
\end{theorem}

In the case of even-dimensional manifolds the following simplified description of $\Spin^c$-structures will be useful for us.  

\begin{definition} 
\label{reduced-spinordef}
Let $M^n$ be a smooth, even-dimensional manifold. A \emph{reduced $\Spin^c$-structure} on $M$ consists of a Riemannian metric on $M$ and a Dirac bundle $S$ ($\Z/2$-graded, but with no $n$-grading structure) whose fiber dimension is $2^{\frac n2}$.  We shall call $S$ a \emph{reduced spinor bundle}.
\end{definition}

If $n$ is even then the complex Clifford algebra $\C_n$ is isomorphic to the matrix algebra $M_{2^{\frac n 2}} (\C)$, and hence has a unique representation $V_n$ of dimension $2^{\frac  n2}$.  The operator 
$$
\gamma = i^{\frac n 2} e_1\cdots e_n
$$
provides $V_n$ with a $\Z/2$-grading.  If $S$ is a reduced spinor bundle, as in the definition, then the tensor product 
$ S \hat \otimes V_n 
$ 
is a spinor bundle in the sense of Definition~\ref{spinordef}, and conversely every spinor bundle in the sense of Definition~\ref{spinordef} is of this form. If we temporarily denote by $[M]_{\operatorname{red}}\in K_0(M)$ the $K$-homology class of the Dirac operator on the reduced spinor bundle $S$, then under the periodicity map $K_0(M)\to K_n(M)$ the fundamental class $[M]_{\operatorname{red}}$ maps to $[M]$.

We conclude this section by comparing reduced spinor bundles  with the Dirac bundles $\Cliff _{\frac 12} (TN)$ that we introduced in Section~\ref{dirac-op-sec}.

Let $N$ be an even-dimensional, oriented Riemannian manifold, and assume it admits a $\Spin^c$-structure, with reduced spinor bundle $S$.  As we noted above, the complex Clifford algebra of a Euclidean vector space of dimension $n=2k$ is isomorphic to the algebra of complex $2^k\times 2^k$ matrices. It follows by counting dimensions that the natural map 
$
\Cliff (TN) \to \End (S)
$
is an isomorphism.  Hence there is an isomorphism
$$
\Cliff (TN) \cong S\hat \otimes S^* 
$$
compatible with the left and right actions by Clifford multiplication.   

\begin{proposition}
\label{half-cliff-calc-prop}
Let $S$ be a reduced spinor bundle for $N$ and denote by $\S^*_{+}$  the even-graded part of its dual.
There is an isomorphism of Dirac bundles 
$$\Cliff_{\frac 12}(TN)\cong S \otimes S^*_+{}.
$$ 
\end{proposition}

\begin{proof}
The reduced spinor bundle determines a full spinor bundle for $M$, which in turn determines the orientation of $M$, as described earlier.  Having fixed this orientation, the operator $\gamma$ acts as $+1$ on $S_{+}$ and $-1$ on $S_{-}$.   So the proposition follows from the isomorphism $\Cliff (TN) \cong S\hat \otimes S^* $.
\end{proof}

\section{Review of Geometric K-Homology}

\begin{definition}
Let $X$ be a paracompact Hausdorff space and let $Y$ be a closed subspace of $X$.  A \emph{$K$-cycle} for the pair $(X,Y)$ is a triple $(M,E,\phi)$ consisting of:
\begin{enumerate}[\rm (i)]
\item  A smooth, compact manifold $M$ (possibly with boundary), equipped with a $\Spin^c$-structure.
\item A smooth, Hermitian vector bundle $E$ on $M$.
\item A continuous map $\phi\colon M \to X$ such that $\phi[ \partial M] \subseteq Y$.
\end{enumerate}
\end{definition}

\begin{remark}
The manifold $M$ need not be connected.  Moreover the components of $M$ may have differing dimensions.
\end{remark}

Two $K$-cycles are \emph{isomorphic} if there are compatible isomorphisms of all of the above three components in the definition of $K$-cycle (this includes an isomorphism of spinor bundles).  Following \cite{baum-douglas82} we are going to construct an abelian group from  sets of isomorphism classes of cycles so as to obtain  ``geometric'' $K$-homology groups for the pair $(X,Y)$.   In order to define the relations in these groups we need to introduce several kinds of operations and  relations involving $K$-cycles.

\begin{definition}
If $(M,E,\phi)$ and $(M',E',\phi')$ are two $K$-cycles for $(X,Y)$, then their \emph{disjoint union} is the $K$-cycle   $(M\cup M', E\cup E', \phi\cup \phi')$.
\end{definition}

\begin{definition} 
If $(M,E,\phi)$ is a $K$-cycle for $(X,Y)$, then its \emph{opposite} is the $K$-cycle $(-M, E, \phi)$, where $-M$ denotes the manifold $M$ equipped with the opposite $\Spin^c$-structure. 
\end{definition} 

\begin{definition}
\label{bordism-def}
A \emph{bordism} of  $K$-cycles for the pair $(X,Y)$ consists of the following data:
\begin{enumerate}[\rm (i)]
\item A smooth, compact manifold $L$, equipped with a $\Spin^c$-structure.
\item A smooth, Hermitian vector bundle $F$ over $L$.
\item A continuous map $\Phi \colon L\to X$.
\item A smooth map $f\colon \partial L \to \R$ for which $\pm 1 $ are regular values, and for which $\Phi [ f ^{-1} [ -1,1]]\subseteq Y$.
\end{enumerate}
\end{definition}

To understand the definition, it is best to consider the case where $Y= \emptyset$.   In this case it follows from condition (iv) that the set $f[-1,1]$ is empty, and therefore the boundary of $L$ is divided by $f$ into two components: $M_{+} = f^{-1} (+1,+\infty)$ and $M_{-} = f^{-1} (-\infty, -1)$. 
We therefore obtain two $K$-cycles $(M_{+}, F|_{M_{+}}, \Phi|_{M_{+}})$ and 
$(M_{-}, F|_{M_{-}}, \Phi|_{M_{-}})$, and we shall say that the first is \emph{bordant} to the opposite of the second.  

In the case where $Y$ is non-empty the sets $M_{+} = f^{-1} [+1,+\infty)$ and $M_{-} = f^{-1} (-\infty, -1]$ are manifolds with boundary, and we obtain, as before two $K$-cycles $(M_{+}, F|_{M_{+}}, \Phi|_{M_{+}})$ and 
$(M_{-}, F|_{M_{-}}, \Phi|_{M_{-}})$, but now for the pair $(X,Y)$.  Once again we shall say that the first is bordant to the opposite of the second.  

The purpose of the function $f$ in Definition~\ref{bordism-def} is to provide a notion of bordism for manifolds with boundary without having to introduce manifolds with corners.  Bordism is an equivalence relation.  

We have one more operation on $K$-cycles to introduce.
Let $M$ be a $\Spin^c$-manifold and let $W$ be a $\Spin^c$-vector bundle over $M$.  Denote by $\mathbf 1$ the trivial, rank-one real vector bundle.    The direct sum $W\oplus \mathbf 1$ is a $\Spin^c$-vector bundle, and moreover the total space of this bundle may be equipped with a $\Spin^c$ structure in a canonical way, up to concordance.  This is because its tangent bundle fits into an exact sequence 
$$
\xymatrix{
0 \ar[r] & \pi^*[W\oplus \mathbf 1] \ar[r] & T(W\oplus \mathbf 1)\ar[r] & \pi ^* [ TM] \ar[r] & 0 , 
}
$$
where $\pi$ is the projection from $W\oplus \mathbf 1$ onto $M$, so that, upon choosing a splitting, (or equivalently, choosing a Riemannian metric on the manifold $W\oplus \mathbf 1$ which is compatible with the above sequence) we have a direct sum decomposition 
$$
 T(W\oplus \mathbf 1) \cong \pi^*[W\oplus \mathbf 1] \oplus  \pi ^* [ TM] .
$$
Different splittings result in concordant $\Spin^c$-structures. 

Let us now denote by $Z$ the unit sphere bundle of the bundle $W\oplus \mathbf 1$.
Since $Z$ is the boundary of the disk bundle, we may equip it with a natural $\Spin^c$-structure by first restricting  the given $\Spin^c$-structure on total space of $W\oplus \mathbf 1$  to the disk bundle, and then taking the boundary of this $\Spin^c$-structure to obtain a $\Spin^c$-structure on the sphere bundle.

\begin{definition}
\label{EF-def}
Let $(M,E,\phi)$ be a $K$-cycle for $(X,Y)$ and let $W$ be a $\Spin^c$-vector bundle over $M$ with even-dimensional fibers.  Let $Z$ be the sphere bundle of $W\oplus \mathbf 1$, as above. The vertical tangent bundle of  $Z$  has a natural $\Spin^c$-structure (one applies the boundary construction of Definition~\ref{boundary-S-def} to the pullback of $W\oplus \mathbf 1$ to $Z$).  Denote by $S_V$ the corresponding reduced spinor bundle and let 
$
F = S_{V, +}^* 
$.
In other words, define $ F$ to be the dual of the even-graded part of the $\Z/2$-graded bundle $S_V$.
  The \emph{modification} of $(M,E,\phi)$ associated to $W$ is the $K$-cycle
$(Z,  F\otimes \pi^*E , \phi\circ \pi)$.
\end{definition}

We are now ready to define the Baum-Douglas geometric $K$-homology groups.

\begin{definition}
\label{k-geom-def}
Denote by $K^{\geom}(X,Y)$ the set of equivalence classes of $K$-cycles over $(X,Y)$, for the equivalence relation generated by the following relations:
\begin{enumerate}[\rm (i)]
\item If $(M,E_1,\phi)$ and $(M,E_2,\phi)$ are two $K$-cycles with the same $\Spin^c$-manifold $M$ and map $\phi\colon M\to X$, then 
$$ (M\cup M ,E_1\cup E_2 ,\phi\cup \phi ) \sim  (M, E_1\oplus E_2, \phi).
$$
\item 
If $(M_1,E_1,\phi_1)$ and $(M_2,E_2,\phi_2)$ are bordant $K$-cycles then 
$$(M_1,E_1,\phi_1) \sim (M_2,E_2,\phi_2).$$
\item If $(M,E,\phi)$ is a $K$-cycle, and if $W$ is an even-dimensional $\Spin^c$-vector bundle over $M$, then 
$$
(M,E,\phi) \sim (Z, F\otimes \pi^*E, \phi\circ \pi),
$$
where $(Z, F\otimes \pi^*E, \phi\circ \pi)$ is the modification of $(M,E,\phi)$ given in Defintion~\ref{EF-def}.
\end{enumerate}
\end{definition}

The set $K^{\geom}(X,Y)$ is in fact an abelian group.  The addition operation is given by disjoint union,
$$
 [ M_1,E_1,\phi_1]  + [M_2,E_2,\phi_2] = [ M_1 \cup M_2, E_1\cup E_2, \phi_1\cup \phi_2],
 $$
 and the additive inverse of a cycle is obtained by reversing the $\Spin^c$-structure: 
 $$
 - [M,E,\phi] = [-M,E,\phi].
 $$
 The neutral element is represented by the empty manifold, or any cycle
 bordant to the empty manifold.
 
\begin{definition}
Denote by $K^{\geom}_{\ev}(X,Y)$ and $K^{\geom}_{\odd}(X,Y)$ the subgroups of the group $K^{\geom}(X,Y)$ composed of equivalence classes of $K$-cycles $(M,E,\phi)$ for which every connected component of $M$ is even dimensional and odd dimensional, respectively.
\end{definition}

The groups $K_{\ev/\odd}^{\geom}(X,Y)$ are functorial in $(X,Y)$, and they satisfy weak excision: if $U$ is an open subset of $Y$ whose closure is in the interior of $Y$, then 
$$
K_{\ev/\odd} ^{\geom}(X\setminus U, Y\setminus U ) \cong K_{\ev/\odd}^{\geom} (X,Y).
$$
There is moreover a ``homology sequence'' 
$$
\xymatrix{
K_{\ev} ^{\geom} (Y) \ar[r] & K_{\ev}  ^{\geom}(X) \ar[r] & K_{\ev}  ^{\geom}( X,Y) \ar[d] \\
K_{\odd} ^{\geom}(X,Y) \ar[u] & K_{\odd}  ^{\geom}(X)\ar[l] & K_{\odd} ^{\geom}(Y)\ar[l]
}
$$
(where as usual we define $K_{\ev}^{\geom}(Y) = K_{\ev}^{\geom} (Y,\emptyset)$, and so on).  The boundary maps take a $K$-cycle $(M,E,\phi)$ for $(X,Y)$ to the boundary cycle $(\partial M, E|_{\partial M}, \phi | _{\partial M})$ for $Y$ (it is easily verified that this definition  is compatible with the equivalence relation used to define the geometric $K$-homology groups).  The composition of any two successive arrows is zero.  However it is not obvious that the sequence is exact.  For the special  
case of finite $CW$-pairs this exactness will follow from the main theorem of the paper, which identifies geometric $K$-homology with Kasparov $K$-homology.

\section{Natural Transformation and Formulation of the Main Theorem}

Now let $(X,Y)$ be a pair of compact and metrizable spaces.  We associate to each $K$-cycle $(M,E,\phi)$ for $(X,Y)$ a class $\langle M,E,\phi\rangle $ in Kasparov $K$-homology, as follows.  Denote by $M^\circ $ the interior of $M$, which is an open $\Spin^c$-manifold.  The $\Spin^c$-structure on $M$ determines a spinor bundle $\S$ on $M^\circ$ by restriction, and of course the complex vector bundle $E$ also restricts to $M^\circ$.  The tensor product $\S \otimes E$ is a Dirac bundle over $M^\circ$, and if $D_E$ is an associated Dirac operator, then we can form the class 
$$
 [ D_E ] \in K_n ( M^\circ)
 $$
 (here $n$ is the dimension of $M$).
The map $\phi\colon M\to X$ restricts to a proper map from $M^\circ$ into $X\setminus Y$, and we can therefore form the class
$$
 \phi_* [ D_E] \in K_n (X,Y).
$$

\begin{theorem}
The correspondence $(M,E,\phi)\mapsto     \phi_* [ D_E] $ determines a functorial map 
$$
\mu\colon K^{\geom}_{\ev/\odd} (X,Y) \to K_{\ev/\odd}(X,Y)
$$
which is compatible with boundary maps in geometric and analytic $K$-homology.
\end{theorem}

\begin{proof}
The only thing to check is that the correspondence is compatible with the relations in Definition~\ref{k-geom-def} which generate the equivalence relation on cycles used to define geometric $K$-homology.  Once this is done, functoriality will be clear from the construction of $\langle M,E,\phi\rangle $ and compatibility with boundary maps will follow from Theorem~\ref{bdiracdirac}. 

Compatibility with relation (i) from Definition~\ref{k-geom-def} is straightforward.  Compatibility with relation (ii) follows from Theorem~\ref{bdiracdirac}.  So the proof reduces to showing that the correspondence is compatible with the relation (iii) of vector bundle modification.

Let $(M,E,\phi)$ be a $K$-cycle for $(X,Y)$ and let $n=\dim (M)$ (by working with one component of $M$ at a time we can assume that $\dim (M)$ is well-defined).  Let $W$ be a $\Spin^c$-vector bundle over $M$ of even fiber dimension $2k$.  Let $\S_M$ be the spinor bundle for $M$, and let $S_V$ be the \emph{reduced} spinor bundle for the vertical tangent bundle of the sphere bundle $\pi\colon Z\to M$.  Form the tensor product
$$
 S_{Z} =  S_V \hat \otimes \pi^*[ \S_M] .
$$
This is neither a fully multigraded spinor bundle for $Z$ nor a reduced spinor bundle, but something in between.  If $D_Z$ is a Dirac operator for $S_Z$ then the class $[D_Z]\in K_n(Z)$ is the image of the $K$-homology fundamental class   $[Z]\in K_{n+2k}(Z)$ under the periodicity isomorphism
$ K_{n+2k}(Z) \cong K_n(Z)$.
Similarly, if $D_{Z,  F\otimes \pi * E}$ is a Dirac operator for the tensor product bundle 
$ 
S_{Z} \otimes F  \otimes \pi^*E $, then the class $[D_{Z,  F\otimes \pi * E}]\in K_n(Z)$ is the image of the $K$-homology class of the modification $(Z, F\otimes \pi^*E, \phi_Z)$ of the cycle $(M,E,\phi)$  under the same periodicity isomorphism.

To prove compatibility with the relation (iii) in the definition of geometric $K$-homology we need to show that $[D_{Z,  F\otimes \pi * E}]$ is equal to the class $[D_{M,E}]\in K_n (M)$.  But writing  
$$
 S_{Z} \otimes  F \otimes \pi^*E\cong  [S_V\otimes F]  \hat \otimes\pi^* [S_M\otimes E] ,$$
we see that this follows from Propositions~\ref{fiber-bundle-prop}, \ref{index-calc-prop} and \ref{half-cliff-calc-prop}.\end{proof}

We can now state the main theorem in this paper.

\begin{theorem}
\label{main-theorem}
If $(X,Y)$ is a finite $CW$-pair then the homomorphism
$$
\mu\colon K^{\geom}_{\ev/\odd} (X,Y) \to K_{\ev/\odd} (X,Y)
$$
is an isomorphism.
\end{theorem}

The proof will be carried out in the remaining sections. 

\section{Outline of the Proof}

We wish to prove that if $X$ is a finite $CW$ complex, then the homomorphisms
$$
\mu \colon K_{\ev/\odd} ^{\geom} (X) \to K_{\ev/\odd} (X)
$$
are isomorphisms.  What makes this tricky is that we don't yet know that geometric $K$-homology is a homology theory.   To get around this problem we are going to define a  ``technical'' homology theory $k_{\ev/\odd}(X,Y)$ which fits into a commuting diagram
$$
\xymatrix{
k_{\ev/\odd}(X,Y) \ar[rr]^{\alpha} \ar[dr]_{\beta}& & K_{\ev/\odd} (X,Y) \\
& K^{\geom}_{\ev/\odd} (X,Y) \ar[ur]_{\mu} &
} 
$$
in which the horizontal arrow is a natural transformation between homology theories.
Having done so, the proof will be completed in two steps: 
\begin{enumerate}[\rm (a)]
\item We shall check that when $X$ is a point and $Y$ is empty, the horizontal arrow is an isomophism.  It will follow that the horizontal arrow is an isomorphism for every finite $CW$ pair $(X,Y)$.
\item We shall prove that for every finite $CW$ pair $(X,Y)$, the map in the diagram from $k_{\ev/\odd} (X,Y) $ to $K_{\ev/\odd}^{\geom} (X,Y)$ is surjective. 
\end{enumerate}
It is clear that (a) and (b) together will imply that all the arrows in the diagram are isomorphisms, for every finite $CW$ pair $(X,Y)$.

The reader who is acquainted with the definition of $K$-homology starting from the Bott spectrum will see that our definition of $k_{\ev/\odd} (X,Y)$ is extremely close to the spectrum definition of $K$-homology.  However the definition which is presented in the next section is not designed with this in mind.

\section{Definition of the Technical Group}

Fix a model $\K$ for the $0$th space of the Bott spectrum,
e.g.~$\mathbb{Z}\times BU$. We shall use the following features of this space:
\begin{enumerate}[\rm (a)]
\item 
If $X$ is a pointed finite $CW$ complex, then there is a natural isomorphism
$$
 K^0(X,*)\cong [ X, \K]^{+}
 $$
between the relative Atiyah-Hirzebruch $K$-theory group $K^0(X,*)$ and the set
of 
homotopy classes of maps from $X$ into $\K$. Here $*$ is the base point of $X$
and $[ X, \K]^{+}$ denotes the
set of homotopy classes of basepoint-preserving  maps (recall that $\K$ is a
base-pointed space). 
\item    
There is a basepoint-preserving  map $m\colon \K\wedge \K \to \K$ which induces the operation of tensor product (the ring structure) on $K^0(X)$.
\end{enumerate}

\begin{example} 
We could take $\K$ to be the space of all Fredholm operators on a separable, infinite-dimensional Hilbert space $H$ (the Fredholm operators are topologized by the operator-norm topology).  The isomorphism (a) is described in \cite{atiyah}.   For later use we note that the set of connected components of $\K$ is isomorphic to $\Z$, the isomorphism being given by the Fredholm index.\end{example}

\begin{definition}
Let $S^2$ be the standard $2$-sphere, equipped with its standard $\Spin^c$ structure as the boundary of the ball in $\R^3$.   Denote by $\beta\colon \S^2 \to \K$ a basepoint-preserving map which, under the isomorphism $[S^2 , \K]^{+} \cong K^0(\R^2)=K^0(S^2,*)$, corresponds to the difference 
$[S_{+}^*] -[\mathbf 1]$.  Here $S^*_{+}$ is the dual of the positive part of the reduced spinor bundle on the $2$-sphere (which is the one-point compactification of $\R^2$), and $\mathbf 1$ is the trivial line bundle.
\end{definition} 

\begin{remark}
Note that $S^*_{+}$ is a line bundle, so that the difference $[S_{+}^*] -[\mathbf 1]$ has virtual dimension zero. 
\end{remark}

Now, we are going to construct the ``technical'' homology groups $k_{\ev/\odd} (X,Y)$ using the space $\K$, the map $\beta$, and the notion of \emph{framed bordism}, which we briefly review. 

\begin{definition} 
A \emph{framed manifold} is a smooth, compact  manifold $M^n$ with a given
 stable trivialization\footnote{In this definition we are using $\mathbf n$ or
 $\mathbf k$ to denote the trivial \emph{real} vector bundle of rank $n$ or
 $k$, respectively. } 
 of its tangent bundle:
$$
 \mathbf k \oplus TM \cong \mathbf {k}\oplus \mathbf{n}
$$
We shall identify two stable trivializations if they are stably homotopic (that is, homotopic after forming the direct sum with the identity map on an additional trivial summand).  Thus a framed manifold is a smooth,  compact manifold together with a stable homotopy class of stable trivializations of its tangent bundle.
\end{definition}

\begin{definition} 
If $(X,Y)$ is any paracompact and Hausdorff pair then we shall denote by $\Omega^F_n(X,Y)$ the $n$-th framed bordism group of the pair $(X,Y)$.  Thus $\Omega_n^F (X,Y)$ is the set of all bordism classes of maps from framed manifolds into $X$, mapping the manifold boundaries into $Y$.  Compare Definition~\ref{bordism-def} or \cite{MR40:2108}.  Note that the boundary of a framed manifold $M$ is itself a framed manifold in a natural way: starting from a stable trivialization
$$
 \mathbf k \oplus TM \cong \mathbf {k}\oplus \mathbf{n}
$$
we use an inward pointing normal vector field on $\partial M$ to obtain a stable trivialization
$$
 \mathbf k \oplus \mathbf 1\oplus T\partial M \cong \mathbf {k}\oplus \mathbf{n}.
$$
We use the inward pointing normal to agree with orientation conventions established earlier.
\end{definition}

We can now define our ``technical'' homology theory $k_{\ev/\odd} (X,Y)$.  For a finite $CW$ pair $(X,Y)$ and an integer $n$, form a direct system of abelian groups 
$$
 \Omega_{n}^F(X\times \K, Y\times \K) \to 
\Omega_{n+2}^F(X\times \K, Y\times \K) \to 
\Omega_{n+4}^F(X\times \K, Y\times \K) \to 
\cdots
$$
as follows. Given a cycle $f\colon M\to X\times \K$ for $\Omega_{n+2k} ^F(X\times \K, Y\times \K)$, the composition
$$
\xymatrix{
M\times S^2 \ar[r]^-{f\times \beta} & X\times \K\times \K \ar[r]^-{1\times m} &X\times \K}
$$
is a cycle for $\Omega_{n+2k+2}^F (X\times \K, Y\times \K)$. This defines the map from $\Omega_{n+2k}^F(X\times \K, Y\times \K)$ to $ \Omega_{n+2k+2}^F(X\times \K, Y\times \K)$ which appears in the directed system. 

\begin{definition}
Denote by $k_{\ev/\odd} (X,Y)$ the direct limit of the above directed system, for $n$ even/odd.
\end{definition}

Since $\Omega^F_*$ is itself a homology theory (on finite $CW$ pairs), and since direct limits preserve exact sequences, it is clear  that $k_*$ is a homology theory. 

The map $\beta$ from $k_{\ev/\odd} (X,Y)$ into $K^{\geom}_{\ev/\odd}(X,Y)$ which appears in Section~7 is defined as follows (for notational simplicity we will only describe the construction for the absolute groups  $k_{\ev/\odd}(X)$, not the relative groups).  If $M$ is a framed manifold then the framing 
$ \bold k \oplus TM \cong \bold k\oplus \bold n$ determines a $\Spin^c$-structure on $M$.  A map $M\to X\times \K$ determines a map $\phi\colon M\to X$ and a $K$-theory class for $M$, which we may represent as a difference $[E_1]- [E_2]$ for some vector bundles $E_1$ and $E_2$.  A map  $$ 
\beta_n \colon  \Omega_{n}^F(X\times \K) \to K^{\geom}_{n}(X)$$ is defined by associating to the bordism class of  $M\to X\times \K$ the difference of $K$-cycles $(M,E_1,\phi)  - (M,E_2,\phi)$.   It follows from part (i) of Definition~\ref{k-geom-def} that the $K$-homology class of this difference does not depend on the choice of $E_1$ and $E_2$ to represent the $K$-theory class on $M$.  It follows from part (ii) of the definition that the $K$-homology class only depends on the bordism class of the map $M\to X\times \K$.  Finally, it follows from part (iii) of the definition that the diagram
$$\xymatrix{\Omega^F_{n} (X\times \K)\ar[d] \ar[r]^-{\beta_n}&K_{n}^{\geom}(X)\ar[d]^{=} \\
\Omega ^F_{n+2} (X\times \K)  \ar[r]_-{\beta_{n+2}} & K_{n+2}^{\geom} (X)}
$$
is commutative (on the right is the periodicity isomorphism  described in Section~2). 
Since $\beta$ is  compatible with  the direct limit procedure using which $k_{\ev/odd}(X)$ is obtained from the framed bordism groups, we obtain  maps 
$$
 \beta  \colon k_{\ev/\odd}(X) \to K^{\geom}_{\ev/\odd}(X)$$
 as required.

\section{Proof of the Main Theorem}

\subsection{Proof of (a)} 
 We wish to show that the maps $
\alpha \colon k_{n}(\Pt) \to K_{n}(\Pt)$ are isomorphisms for $n=0$ and $n=1$.  If $W$ is any (base-pointed) topological space then by the Pontrjagin-Thom isomorphism \cite{MR4:249e,MR22:5980,MR15:890a} the $n$th framed bordism group of $W$ is isomorphic to the $n$th stable homotopy group of $W$: 
$\Omega_n^F (W) \cong \pi_n^S(W)
$.
According to Bott periodicity, \cite{MR21:1588,MR31:2727} the second loop space of $\K$ has the homotopy type of $\K$.  In fact the map
$$
 \xymatrix{ S^2\wedge \K \ar[r]^{\beta\wedge 1} & \K\wedge \K\ar[r]^{m}& \K}$$
 induces a homotopy equivalence $\K\sim \Omega^2 \K$.  This, and the fact that the flip map $S^2\wedge S^2\to S^2\wedge S^2$ is homotopic to the identity map, imply that the evident maps 
 $$\varinjlim \pi_{n+2k} (\K)\to \varinjlim \pi_{n+2k}^S(\K)$$
are  isomorphisms. The first direct limit is formed by associating to a map $f\colon S^{n+2k}\to \K$ the composition 
$$\xymatrix{
S^2\wedge S^{n+2k}\ar[r]^{1\wedge f} & S^2\wedge \K \ar[r]^{m}& \K,}$$
and the second direct limit is formed using  a similar procedure, starting with maps from $S^{n+2k+2j} $ into $S^{2j}\wedge \K$.   To verify the assertion, view the second direct limit as the limit of the array
$$
\xymatrix{
\vdots &  \vdots & \vdots&\\
\pi_6 (S^4\wedge \K) \ar[u]^s\ar[r]^b &\pi_8 (S^4\wedge \K) \ar[u]^s\ar[r]^b &\pi_{10} (S^4\wedge \K) \ar[u]^s\ar[r]^b & \dots\\
\pi_4 (S^2\wedge \K) \ar[u]^s\ar[r]^b &\pi_6 (S^2\wedge \K) \ar[u]^s\ar[r]^b &\pi_8 (S^2\wedge \K) \ar[u]^s\ar[r]^b & \dots\\
\pi_2 ( \K) \ar[u]^s\ar[r]^b &\pi_4 ( \K) \ar[u]^s\ar[r]^b &\pi_6 ( \K) \ar[u]^s\ar[r]^b & \dots\\
}
$$
in which the vertical maps are suspension by $S^2$ and the horizontal maps are induced from suspension by $S^2$, followed by composition with  
  $b\colon S^2 \wedge \K \to \K$ defined by
$$
\xymatrix{
S^2\wedge \K \ar[r]^{\beta\wedge 1}& \K\wedge \K \ar[r]^{m} &\K}.
$$
The first direct limit is then the direct limit of the bottom row, and the required isomorphism follows from these facts:
\begin{enumerate}[\rm (i)]
\item
If $x\in \pi_{2 k}(S^{2j}\wedge \K)$, and if $s(x)=0$, then $b(x)=0$.
\item 
If $x\in \pi_{2 k}(S^{2j}\wedge \K)$ for some $j>0$, and if $x=b(y)$, for some $y\in \pi_{2k-2} (S^{2j}\wedge \K)$, then $x= s(z)$, for some $z\in \pi_{2k-2}( S^{2j-2}\wedge \K)$.
\end{enumerate}
Item (i) is an immediate consequence of the definition of the map $b$.  As for item (ii), if $x = b(y)$, then $x$ can be written as a composition
$$
\xymatrix{
S^2\wedge S^{2k-2} \ar[r]^-{1 \wedge y}& S^2 \wedge S^{2j}\wedge \K \ar[r]^-{1\wedge b} & S^{2j}\wedge \K}.
$$
Writing $S^{2j}$ as $S^{2j-2} \wedge S^{2}$, and using the fact that the flip on $S^2 \wedge S^2$ is homotopic to the identity, we can write this composition as 
$$
\xymatrix{
S^2\wedge S^{2k-2} \ar[r]^-{1\wedge y} & S^2 \wedge S^{2j-2}\wedge S^2\wedge \K  \ar[r]^-{1\wedge 1 \wedge b} & S^2\wedge S^{2j-2}\wedge \K}.
$$
This is clearly in the image of the map $s$.

Now 
 $$
 k_n (\Pt) = \varinjlim \Omega^F_{n+2k} (\K) \cong \varinjlim \pi^S_{n+2k}(\K)\cong \varinjlim \pi_{n+2k} (\K)\cong \pi_n (\K).$$
  As a result we obtain the isomorphisms
$$
 k_n (\Pt)  \cong \pi_{n } (\K)\cong K^0(\R^n)$$
which implies that $k_{\ev}(\Pt)  \cong \Z$ and $k_{\odd} (\Pt) = 0$.   It follows immediately that the map $\alpha\colon k_{\odd} (\Pt) \to K_{\odd} (\Pt)$ is an isomorphism, since both domain and range are zero.  In the even case  the map 
$$
k_0 (\Pt) \cong  \pi_{0} (\K) \to K_0(\Pt) \cong \Z
$$
sends a Fredholm operator $T$ to the index of $T$.  This map is an isomorphism.

\subsection{Proof of (b)}
We wish to prove that the map $k_{\ev/\odd} (X,Y) \to K^{\geom}_{\ev/\odd}(X,Y)$ is surjective.  The image of this map consists precisely of the equivalence classes of $K$-cycles $(N,F,\psi)$ for which $N$ is a framed $\Spin^c$-manifold.   So we must prove that if  $(M,E,\phi)$ is any  $K$-cycle for $(X,Y)$, then there is an equivalent $K$-cycle $(N,F,\psi)$  for which $N$ is a framed $\Spin^c$-manifold. 

To do this, choose a smooth real vector bundle $V$, with even-dimensional fibers, such that $TM\oplus V $ is trivializable, and fix an isomorphism 
$$TM \oplus V \cong \mathbf n\oplus \mathbf k. $$
The trivial bundle $\mathbf n\oplus \mathbf k$ has a canonical $\Spin^c$-structure, and the above isomorphism and the following lemma therefore define a $\Spin^c$-structure on $V$.  

\begin{lemma} 
Let $V$ and $W$ be real, orthogonal vector bundles over the same space $X$.  Assume that $V$ and $V\oplus W$ are equipped with $\Spin^c$-structures.  There is a $\Spin^c$-structure on $W$ whose direct sum with the given $\Spin^c$-structure on $V$ is the given $\Spin^c$-structure on $V\oplus W$.
\end{lemma} 

\begin{proof} 
Let $S_V$ be a (non-reduced) spinor bundle for $V$ and let $S_{V\oplus W}$ be
the same for $V\oplus W$.   Denote by $\S_{W}$ the bundle of fiberwise linear
maps $S_V\to S_{V\oplus W}$ which graded-commute with the Clifford action of
$V$ and which graded-commute with the action of the first $k$ multigrading
operators $\varepsilon_1, \dots ,\varepsilon_k$, where $k= \rank (V)$.  The
bundle $W$ acts on $S_{V\oplus W}$, as do the remaining multigrading
operators. By composition, $W$ and the remaining multigrading operators also
act on $S_W$.  A local consideration shows that we obtain a (non-reduced)
spinor bundle for $W$, and that it has the required property with respect to
direct sum.
\end{proof}

The vector bundle modification of the $K$-cycle $(M,E,\phi)$ by $V$ is a $K$-cycle whose manifold is framed, as required.
  
\section{Appendix: The Real Case} 

In this appendix we briefly discuss the changes needed to prove the result analogous to Theorem~\ref{main-theorem} in $KO$-homology.  

Kasparov's theory readily adapts to the real case.  A real Hilbert space can be viewed as a complex Hilbert space equipped with a conjugate-linear isometric involution.  A real $C^*$-algebra is the same thing as a complex $C^*$-algebra equipped with a conjugate linear involutive $*$-automorphism (which, unlike the $*$-operation, preserves the order of products).  By including these complex-conjugation operators, the definitions of Section~\ref{kasp-sec} extend immediately to the real case.   The only difference is that in the real case the counterpart of the formal periodicity map $K^{-p} \to K^{-(p+2)}$ does not exist.  However the four-fold composition of this map is compatible with real structures, and defines a real formal periodicity map 
$KO^{-p} \to KO^{-(p+8)}$.  The results of Section~\ref{dirac-op-sec} carry over without change, except that the bundle $\Cliff_{\frac 12} (TN)$ is a real Dirac bundle only when the dimension of $N$ is a multiple of $4$.   Our discussion of $\Spin^c$ structures in Section~\ref{spinc-sec}  is designed to carry over to the real case just by replacing complex Clifford algebras with real Clifford algebras; reduced real spinor bundles exist in dimensions which are multiples of $8$.  The geometric definition of $K$-homology is based on $\Spin$-manifolds---the real counterparts of  $\Spin^c$-manifolds---and the real counterpart of Theorem~\ref{main-theorem} is now easy to formulate.   The only really new aspect of the proof is that a more careful treatment of part (a) is required.  The argument given above shows that 
$$
ko_0(\Pt) \cong \pi_0 (\mathbb K \mathbb O)\cong KO^0(\R^n).
$$
Under these isomorphisms, the map $ko_0(\Pt) \to KO^n(\Pt)$ corresponds to the map $KO^0(\R^n) \to KO_n(\Pt)$ which takes a $K$-theory class $x$ to the index of the Dirac operator on $\R^n$ twisted by $x$.  The fact that this map is an isomorphism is another formulation of Bott Periodicity (compare \cite{MR0228000}).

\bibliographystyle{plain}  
\bibliography{references}       

\noindent P.B.:  Department of Mathematics,  Penn State University, University Park, PA 16802. Email: {\tt baum@math.psu.edu}.

\noindent N.H.:  Department of Mathematics,  Penn State University, University Park, PA 16802. Email: {\tt higson@math.psu.edu}.

\noindent T.S.:  Mathematisches Institut, Georg-August-Universit\"at G\"ottingen,
Bunsenstr.\ 3, 
D-37073 G\"ottingen, Germany. Email:  {\tt schick@uni-math.gwdg.de}.

\end{document}